\documentclass[a4paper,12pt]{amsart}
\usepackage{amssymb}
\usepackage{ifthen}
\usepackage{graphicx}
\usepackage{color,xcolor}
\nonstopmode
\numberwithin{equation}{section}
\newtheorem{thm}{Theorem}
\newtheorem{cor}{Corollary}
\newtheorem{lem}{Lemma}
\newtheorem{prop}{Proposition}

\newtheorem{conj}{Conjecture}
\newtheorem{prob}{Problem}

\theoremstyle{definition}
\newtheorem{defn}{Definition}
\newtheorem{ca}{Case}

\newtheorem{rem}{Remark}

\newenvironment{pf}[1][]{%
 \vskip 1mm
 \noindent
 \ifthenelse{\equal{#1}{}}%
  {{\slshape Proof. }}%
  {{\slshape #1.} }%
 }%
{\qed\medskip}

\newcounter{alphabet}

\newenvironment{Thm}[1][]{\refstepcounter{alphabet}%
\bigskip%
\noindent%
{\bf Theorem \Alph{alphabet}}%
\ifthenelse{\equal{#1}{}}{}{ (#1)}%
{\bf .} \itshape}{\vskip 8pt}

\newenvironment{Lem}[1][]{\refstepcounter{alphabet}%
\bigskip%
\noindent%
{\bf Lemma \Alph{alphabet}}%
{\bf .} \itshape}{\vskip 8pt}

\newcounter{alphabet2}

\newcommand{\IC}{{\mathbb C}}
\newcommand{\ID}{{\mathbb D}}

\def\be{\begin{equation}}
\def\ee{\end{equation}}

\newcommand{\ben}{\begin{enumerate}}
\newcommand{\een}{\end{enumerate}}

\newcommand{\blem}{\begin{lem}}
\newcommand{\elem}{\end{lem}}
\newcommand{\bthm}{\begin{thm}}
\newcommand{\ethm}{\end{thm}}
\newcommand{\bcor}{\begin{cor}}
\newcommand{\ecor}{\end{cor}}
\newcommand{\beg}{\begin{exam}}
\newcommand{\eeg}{\end{exam}}
\newcommand{\begs}{\begin{examples}}
\newcommand{\eegs}{\end{examples}}
\newcommand{\bdefe}{\begin{defn}}
\newcommand{\edefe}{\end{defn}}
\newcommand{\bprob}{\begin{prob}}
\newcommand{\eprob}{\end{prob}}
\newcommand{\bques}{\begin{ques}}
\newcommand{\eques}{\end{ques}}
\newcommand{\bei}{\begin{itemize}}
\newcommand{\eei}{\end{itemize}}
\newcommand{\bcon}{\begin{conj}}
\newcommand{\econ}{\end{conj}}
\newcommand{\bop}{\begin{op}}
\newcommand{\eop}{\end{op}}

\newcommand{\bas}{\begin{assertion}}
\newcommand{\eas}{\end{assertion}}

\newcommand{\bfa}{\begin{fact}}
\newcommand{\efa}{\end{fact}}

\newcommand{\bca}{\begin{ca}}
\newcommand{\eca}{\end{ca}}

\newcommand{\bst}{\begin{step}}
\newcommand{\est}{\end{step}}

\newcommand{\bsca}{\begin{sca}}
\newcommand{\esca}{\end{sca}}

\newcommand{\bcl}{\begin{cl}}
\newcommand{\ecl}{\end{cl}}

\newcommand{\bmlem}{\begin{mlem}}
\newcommand{\emlem}{\end{mlem}}

\newcommand{\bscl}{\begin{scl}}
\newcommand{\escl}{\end{scl}}

\newcommand{\bcons}{\begin{conjs}}
\newcommand{\econs}{\end{conjs}}

\newcommand{\bprop}{\begin{prop}}
\newcommand{\eprop}{\end{prop}}

\newcommand{\br}{\begin{rem}}
\newcommand{\er}{\end{rem}}
\newcommand{\brs}{\begin{rems}}
\newcommand{\ers}{\end{rems}}
\newcommand{\bo}{\begin{obser}}
\newcommand{\eo}{\end{obser}}
\newcommand{\bos}{\begin{obsers}}
\newcommand{\eos}{\end{obsers}}
\newcommand{\bpf}{\begin{pf}}
\newcommand{\epf}{\end{pf}}
\newcommand{\ba}{\begin{array}}
\newcommand{\ea}{\end{array}}
\newcommand{\beq}{\begin{eqnarray}}
\newcommand{\beqq}{\begin{eqnarray*}}
\newcommand{\eeq}{\end{eqnarray}}
\newcommand{\eeqq}{\end{eqnarray*}}

\newcounter{minutes}
\divide\time by 60
\newcounter{hours}
\multiply\time by 60 \addtocounter{minutes}{-\time}

\begin{document}

\bibliographystyle{amsplain}
\title [Bloch and Landau type theorems for pluriharmonic mappings]
{Bloch and Landau type theorems for pluriharmonic mappings} 

\def\thefootnote{}
\footnotetext{ \texttt{\tiny File:~\jobname .tex,
          printed: \number\day-\number\month-\number\year,
          \thehours.\ifnum\theminutes<10{0}\fi\theminutes}
} \makeatletter\def\thefootnote{\@arabic\c@footnote}\makeatother

\author{Ming-Sheng Liu 
}
\address{M.-S. Liu, School of Mathematical Sciences, South China Normal University, Guangzhou, Guangdong 510631, China.}
\email{liumsh65@163.com}

\author{Saminathan Ponnusamy}
\address{S. Ponnusamy, Department of Mathematics,
Indian Institute of Technology Madras, Chennai-600 036, India. }
\address{Department of Mathematics, Petrozavodsk State University, ul., Lenina 33, 185910 Petrozavodsk,
Russia}
\email{samy@iitm.ac.in}

\subjclass[2010]{Primary 31C10; Secondary 32A18, 31B05, 30C65}
\keywords{holomorphic mapping, harmonic mappings, pluriharmonic mappings, Landau-type theorems, Bloch-type theorems
}

\begin{abstract}
In this paper, we establish two new versions of Landau-type theorems for pluriharmonic mappings with  a bounded distortion. 
Then using these results, we derive three Bloch-type theorems of pluriharmonic mappings, which improve the corresponding
results of Chen and Gauthier.
\end{abstract}

\maketitle

\pagestyle{myheadings}
\markboth{M.-S. Liu and S. Ponnusamy}{The Bloch-Landau type theorems for pluriharmonic mappings}

\section{Preliminaries and some basic questions}\label{HLP-sec1}

Let $\mathbb{C}^n$ denote the $n$-dimensional complex Euclidean space so that $\mathbb{C}:=\mathbb{C}^1$, the complex plane.
The conjugate $\overline{z}$ of $z=(z_1,\ldots,z_{n})\in\mathbb{C}^n$ is defined by $\overline{z}=(\overline{z_1},\ldots,\overline{z_{n}})$.
For $z$ and $w=(w_1,\ldots,w_{n})\in\mathbb{C}^n$, we define
\begin{equation*}
\langle{z,w}\rangle:=z_1\overline{w}_1+\cdots+z_{n}\overline{w}_{n} ~\mbox{ and }~ |z|:=\sqrt{\langle{z,z}\rangle}=\left (\sum\limits_{k=1}^n|z_k|^2\right )^{1/2}.
\end{equation*}
For $a\in \mathbb{C}^{n}$,  let $B^{n}(a, r)=\{z=(z_1, \ldots , z_n)\in\mathbb{C}^{n}:~|z-a|<r\}$ be the ball in $\mathbb{C}^{n}$ of radius $r$ with center $a$.
In the case of $n=1$, we use the standard notation $\mathbb{D}(a, r):=B^{1}(a, r)$ so that $\mathbb{D}:=\mathbb{D}(0, 1)$, the open unit disk in $\IC$.
We also let $B^{n}$ denote the unit ball $B^{n}(0,1)$. Evidently, $\mathbb{D}=B^1$ (see \cite{R1980}).

\subsection{Planar harmonic mappings}
For a continuously differentiable complex-valued mapping $f(z)=u(z)+i v(z)$, $z=x+iy$, we use the common notation for its formal derivatives:
\begin{equation}
f_{z}=\frac{1}{2}(f_{x}-if_{y})~\mbox{ and }~ f_{\overline{z}}=\frac{1}{2}(f_{x}+if_{y}).
\end{equation}
We say that $f$ is a harmonic mapping in a simply connected domain $D$ if $f$ is twice continuously differentiable
and satisfies the Laplace equation $\Delta f =4f_{z\,\overline{z}}=0$ in $D$.

Let ${\mathcal H}(\ID)$ denote the set of all harmonic mappings in $\ID$. It is well-known that such  mappings
have the representation $f=h+\overline{g}$, where $h$ and $g$ are analytic functions in $\ID$.  It is convenient to introduce the following   notations:
\beqq
{\mathcal H}_0(\ID)&=& \{f=h+\overline{g}\in {\mathcal H}(\ID):\, g(0)=0\},\\
{\mathcal A}(\ID)&=& \{f=h+\overline{g}\in {\mathcal H}(\ID):\, \mbox{$g(z)\equiv 0$}, h'(0)\neq 0\},\\
{\mathcal A}_0(\ID)&=& \{f\in {\mathcal A}(\ID):\, \mbox{$f(0)=0$}\}
\eeqq

Lewy's theorem \cite{H1936} from 1936 states that a harmonic mapping $f=h+\overline{g}$ is locally univalent on $D$
if and only if the determinant $|J_f|$ of its Jacobian matrix $J_f$ does not vanish on $D$,
where $|J_f|=|f_z|^2-|f_{\overline{z}}|^2=|h'|^2-|g'|^2$. Such a result does not hold in higher   dimensions  (for details see \cite{J1991}).

For a continuously differentiable mapping $f$ on $\ID$,  the maximum and minimum length distortions
of the mapping $f$ are defined respectively by
\begin{eqnarray*}
\Lambda_f (z) &=& \max_{0\leq t\leq2\pi}|f_z(z)+e^{-2it}f_{\overline{z}}(z)|
=|f_z(z)|+|f_{\overline{z}}(z)|, ~\mbox{ and }
\\
\lambda_f (z)&=& \min_{0\leq t\leq2\pi}|f_z(z)+e^{-2it}f_{\overline{z}}(z)|=\big ||f_z(z)|-|f_{\overline{z}}(z)| \big |.
\end{eqnarray*}

Methods of Harmonic mappings have been used to study and solve fluid flow problems (see \cite{AC2012,CM2017}).
For example, in 2012, Aleman and Constantin \cite{AC2012} established a connection between harmonic mappings and ideal fluid flows.
In fact, they have developed ingenious technique to solve the incompressible two dimensional Euler equations in terms of
univalent harmonic mappings. More precisely, the problem of finding all solutions which in Lagrangian variables describing
the particle paths of the flow present a labelling by harmonic mappings is reduced to solve an explicit nonlinear differential
system in $\mathbb{C}^n$ (cf. \cite{CM2017}).


Our primary interest in this paper is to establish several new versions of Landau-type theorems and two improved Bloch-type theorems of pluriharmonic mappings.

\subsection{Landau-Bloch type theorems of harmonic mappings}

The Bloch theorem (1925)  asserts the existence of a positive number $b$ such that for each $f\in {\mathcal A}(\ID)$
there is a disk of radius $b|f'(0)|$ which is the univalent image under $f$ of some subdomain of $\mathbb{D}$. Such a disk is called ``schlicht disk" for $f$.
The supremum of all such numbers $b$ is called the Bloch constant ${\bf B}$.

Bloch's theorem implies the existence of another positive number $\ell$ such that, for each $f\in {\mathcal A}(\ID)$,
$f(\ID)$ contains a disk of radius $\ell |f'(0)|$. The largest possible value of $\ell$, denoted by ${\bf L}$, is known as the Landau constant.
Clearly ${\bf L}\geq {\bf B}$. The exact values of ${\bf B}$ and ${\bf L}$ are not known although the lower and upper bounds are available in the literature.

The classical Landau theorem asserts that if $f\in {\mathcal A}_0(\ID)$ such that $f^{\prime}(0)=1$ and $|f(z)|<M$ for $z\in\mathbb{D}$, then $f$ is univalent
 $\mathbb{D}_{r_0}$, and $f(\mathbb{D}_{r_0})$ contains a disk $\mathbb{D}_{\sigma_0}$, where
$$r_0=\frac{1}{M+\sqrt{M^2-1}} ~\mbox{ and }~ \sigma_0=M r_0^2.
$$
This result is sharp, with the extremal function $f_0(z)=M z \left (\frac{1-M z}{M-z}\right )$.

\bdefe
{\it
A function $f\in {\mathcal H}(\ID)$ is said to belong to ${\mathcal S}_{H}(r;R)$ if  it is univalent in $\ID_r$ and the range
$f(\mathbb{D}_r)$ contains a univalent disk  $\mathbb{D}_R$.
}
\edefe

In 2000, under a suitable restriction,  Chen et al. \cite{CHH2000} established two versions of Landau-type theorems for bounded harmonic
mapping on the unit disk which we now recall them using our notation.

\begin{Thm}\label{Theo-AA} {\rm (\cite[Theorem 3]{CHH2000})}  
 Let $f\in {\mathcal H}_0(\ID)$  such that $f_{\overline{z}}(0)=0$, $f_z(0)=1$, and $|f(z)|<M$ for
$z\in\mathbb{D}$. Then $f\in {\mathcal S}_{H}(r_1;r_1/2)$ with
$$r_1=\frac{\pi^2}{16 m M}\approx\frac{1}{11.105M},
$$
where $m\approx 6.85$ is the minimum of the function $(3-r^2)/(r(1-r^2))$ for $0<r<1$.
\end{Thm}

\begin{Thm}\label{Theo-AB} {\rm (\cite[Theorem 4]{CHH2000})}
Let $f\in {\mathcal H}_0(\ID)$  such that $\lambda_f(0)=1$, and $\Lambda_f(z)\leq\Lambda$ for $z\in\mathbb{D}$. Then $f\in {\mathcal S}_{H}(r_2;r_2/2)$, where
$$r_2=\frac{\pi}{4(1+\Lambda )}.
$$
\end{Thm}

Theorems~A and B 
are not sharp. Better estimates were given in \cite{CHH2011,cpw2011C,CPR2014,CPW2011B,G2006,Huang2014,L2009S,L-CMA2009,LC2018,LLL2019}.
In particular, the sharp version of Theorem \Ref{Theo-AB} were obtained in \cite{Huang2014,L2009S,LC2018}.
For example, the sharp version of Theorem~A 
for $M=1$ reads as follows: 




\vspace{0.2cm}

\begin{Thm}\label{Theo-AD}{\rm (\cite[Theorems 2.4 and 2.5, Remark 2.6]{L-CMA2009})}
Let $f\in {\mathcal H}_0(\ID)$ such that either $J_f(0)=1$ or $\lambda_f(0)=1$, and $|f(z)|<1$ for
$z\in\mathbb{D}$. Then $f\in {\mathcal S}_{H}(1;1)$. The result is sharp.
\end{Thm}

Recently, Liu \cite{LLL2019} proposed the following conjecture which may be regarded as the sharp form of Theorem~A 
for the case $M>1$. 

\bcon\label{Conj-A} {\rm (\cite[Conjecture 3.4]{LLL2019})} 
Let $f\in {\mathcal H}_0(\ID)$  such that $f(0)=0$, $\lambda_f(0)=1$ and $|f(z)|<M$ for $z\in \mathbb{D} $ and for some  $M>1$.
Then $f\in {\mathcal S}_{H}(r_0;\sigma_0)$. The two radii $r_0$ and $\sigma_0$ are sharp, with the extremal mappings
$e^{i\alpha}f_0(z)$ or $e^{i\alpha}\overline{f_0(z)}$, where $\alpha\in\mathbb{R}$ and $f_0(z)=M z \left (\frac{1-M z}{M-z}\right )$.
\econ

\subsection{The Landau-Bloch type theorems of pluriharmonic mappings}
A continuous complex-valued function $\phi$ defined on a domain $\Omega\subset\mathbb{C}^n$ is called a pluriharmonic mapping if, for each fixed $z'\in\Omega$ and $\theta\in\partial B^n$, the function $\phi (z'+\theta\zeta)$ is harmonic in the complex variable $\zeta$, for $|\zeta|$ smaller than the distance
from $z'$ to $\partial B^n$. A mapping $f$ of $\Omega$ into $\mathbb{C}^n$ is called a pluriharmonic mapping if every component of $f$ is pluriharmonic.

A mapping $f$ of $B^n$ into $\mathbb{C}^n$ is pluriharmonic if and only if $f$ has a representation $f = g + \overline{h}$, where $g$ and $h$ are holomorphic mappings (see [4]).

For a continuously differentiable mapping $f:\, B^n\rightarrow \mathbb{C}^m$, $w=f(z)=(f_1(z),\ldots, f_m(z))$, $z=(z_1, \ldots , z_n)$, we denote by $f_z$ and $f_{\overline{z}}$ for the matrices $(\partial f_j/\partial z_k)_{m\times n}$ and $(\partial f_j/\partial \overline{z}_k)_{m\times n}$, respectively.

Denote the   maximum length distortion $\Lambda_f$ and  minimum length distortion $\lambda_f$ by 
\begin{equation*}
\Lambda_f(z) = \max_{\theta \in \partial B^{n}}|f_z(z)\theta+ f_{\overline{z}}(z)\overline{\theta}| ~\mbox{ and }~ \lambda_f(z) = \min_{\theta \in \partial B^{n}}|f_z(z)\theta+ f_{\overline{z}}(z)\overline{\theta}|,
\end{equation*}
respectively, where $\theta$ is regarded as a column vector.
\vspace{0.2cm}

For a continuously differentiable mapping $w = f(z) = (f_1(z), \ldots , f_n(z)), z =(z_1, \ldots , z_n)$, of a domain $\Omega\subset \mathbb{C}^n$ into $\mathbb{C}^n$, let $z_k = x_k + iy_k$ and $f_j = u_j + iv_j$. Note also that $f$ can be regarded as a mapping of a domain in $\mathbb{R}^{2n}$ into $\mathbb{R}^{2n}$. We denote the real Jacobian matrix of this mapping by $J_f$:
\begin{equation*}
J_f=
\begin{pmatrix} \frac{\partial u_1}{\partial x_1}&\frac{\partial u_1}{\partial y_1}&\cdots\cdots&\frac{\partial u_1}{\partial x_n}&\frac{\partial u_1}{\partial y_n}
\\\frac{\partial v_1}{\partial x_1}&\frac{\partial v_1}{\partial y_1}&\cdots\cdots&\frac{\partial v_1}{\partial x_n}&\frac{\partial v_1}{\partial y_n}
\\\vdots&\vdots&\vdots&\vdots&\vdots
\\\frac{\partial u_n}{\partial x_1}&\frac{\partial u_n}{\partial y_1}&\cdots\cdots&\frac{\partial u_n}{\partial x_n}&\frac{\partial u_n}{\partial y_n}
\\\frac{\partial v_n}{\partial x_1}&\frac{\partial v_n}{\partial y_1}&\cdots\cdots&\frac{\partial v_n}{\partial x_n}&\frac{\partial v_n}{\partial y_n}
\end{pmatrix}.
\end{equation*}

Let $B^n$ (resp. $\mathbb{B}^{2n}$) denote the unit ball in $\mathbb{C}^n$ (resp.  $\mathbb{R}^{2n}$). With this notation,
the   maximum length distortion $\Lambda_f$ and minimum length distortion $\lambda_f$ have another equivalent representation:
$$
\Lambda_f=\max\limits_{\theta\in\partial B^{n}}|f_z \theta+f_{\overline{z}} \overline{\theta}|=\max\limits_{\theta\in\partial\mathbb{B}^{2n}}|J_f \theta| ~\mbox{ and }~
\lambda_f=\min\limits_{\theta\in\partial B^{n}}|f_z \theta+f_{\overline{z}} \overline{\theta}|=\min\limits_{\theta\in\partial\mathbb{B}^{2n}}|J_f \theta|.
$$

For a complex or real $n \times n$ matrix $A$, the operator norm of $A$ is defined by
\begin{equation*}
|A| = \sup\limits_{z\neq0}\frac{|Az|}{|z|}=\max_{\theta\in \partial \mathbb{B}^{2n}}|A\theta|.
\end{equation*}
 Thus the maximum distortion $\Lambda_f$ is the $L^2$ operator norm of the Jacobi matrix $J_f$. From now onwards, we identify a point in $\mathbb{C}^n$ or $\mathbb{R}^{2n}$ (real space of dimension $2n$) with a complex or real column vector.

A pluriharmonic mapping $f$ of $B^n$ into $\mathbb{C}^{n}$ is said to be a $K$-mapping if
$$
|J_f (z)| \leq K|{\rm det} J_f (z)|^{1/(2n)}\mbox{ for } z\in B^n.
$$

In 2011, Chen and Gauthier \cite{CHH2011} obtained  Landau theorems and Bloch theorems for pluriharmonic mappings $f: B^n\rightarrow \mathbb{C}^n$.
Recently, Xu and Liu obtained a new version of Landau theorem and a Bloch theorem for pluriharmonic mappings in \cite{XL2020}. Before we recall
the work of Chen and Gauthier, it is convenient to use the following notations.

\bdefe
Define
\beqq
{\mathcal{PH}}(B^{n})&=&\{f:\,B^n \to \mathbb{C}^{n}:\, \mbox{$f$ is pluriharmonic mapping}  \},\\
{\mathcal{PH}}^{\alpha}(B^{n})&=&\{f\in {\mathcal{PH}}(B^{n}):\, {\rm det} J_f (0) = \alpha  \},\\
{\mathcal{PH}}_K(B^{n})&=&\{f:\,B^n \to \mathbb{C}^{n}:\, \mbox{$f$ is pluriharmonic $K$-mapping} \},\\
{\mathcal{PH}}_K^{\alpha}(B^{n})&=&\{f\in {\mathcal{PH}}_K(B^{n}):\, {\rm det} J_f (0) = \alpha  \}, \mbox{ and }\\
{\mathcal{PH}}_{loc,\,K}^{\alpha}(B^{n})&=&\{f\in {\mathcal{PH}}_K^{\alpha}(B^{n}):\, \mbox{$f$ is locally univalent} \}.
\eeqq
A mapping $f\in {\mathcal{PH}}(B^{n})$ is said to belong to ${\mathcal S}_{PH}(r;R)$ if  $f$ is univalent on the ball $B^n (0, r)$ and the range
$f(B^n (0, r))$ covers the ball $B^n(0, R)$.
\edefe

\begin{Thm}\label{Theo-A} {\rm (\cite[Theorem 5]{CHH2011})} 
Let $f\in {\mathcal{PH}}^{\alpha}(B^{n})$ such that $f(0) = 0$, $|f(z)|<M$ for $z\in B^{n}$ and $(4M /\pi)^{2n} \geq \alpha > 0$. Then
$f\in {\mathcal S}_{PH}(\rho_0;R_0)$, where
$$\rho_0 = \frac{\alpha \pi^{2n+1}}{4m(4M)^{2n}} ~\mbox{ and }~ R_0 =\frac{\alpha^2\pi^{4n}}{8m(4M)^{4n-1}},
$$
where $m\approx 4.2$ is the minimum of the function $\frac{2-r^2}{r(1-r^2)}$ on the interval $(0, 1)$.
\end{Thm}

\begin{Thm}\label{Theo-B} {\rm (\cite[Theorem 6]{CHH2011})} 
Let $f\in {\mathcal{PH}}_K^{1}(B^{n})$ with $n > 1$. Then $f(B^{n})$ contains a disk of radius $b_f$, with
\begin{eqnarray*}
b_f\geq R_n' &=& \frac{k_n\pi}{8m}\Big( \frac{k_n\pi}{8K\log (1/(1-k_n))}\Big)^{4n-1}\\
&=&\frac{\pi}{8me(2n-1)}\Big(\frac{\pi}{4}\Big)^{4n-1}\frac{1}{K^{4n-1}}(1+o(1)),
\end{eqnarray*}
where $0<k_n<1$ is the unique number such
$$
4n\log\frac{1}{1-k_n}=(4n-1)\cdot\frac{k_n}{1-k_n}.
$$
\end{Thm}

\begin{Thm}\label{Theo-C} {\rm (\cite[Theorem 7]{CHH2011})} 
Let $f\in {\mathcal{PH}}_{loc,\,K}^{\alpha}(B^{n})$ such that $n > 1$, $(4M /\pi)^{2n} \geq \alpha > 0$, $f(0) = 0$ and $|f(z)|<M$ for $z\in B^{n}$.
Then there exists a domain $\Omega\subset B^n(0, \rho_0')$ such that $0\in\Omega$
and $f$ maps $\Omega$ onto a ball $B^n(0,R_0')$ injectively, where
$$
\rho_0' = \frac{\pi^2 (\alpha)^{1/(2n)}}{16 m M K^{2n-1}},\quad R_0' = \frac{\pi^2 \alpha^{1/n}}{32 m M K^{4n-2}},
$$
and $m$ is the same number as in Theorem~D 
\end{Thm}

\begin{Thm}\label{Theo-D} {\rm (\cite[Theorem 8]{CHH2011})} 
Let $f\in {\mathcal{PH}}_{loc,\,K}^{1}(B^{n})$ and $n > 1$. Then $f(B^{n})$ contains a disk of radius $b_f$, with
$$b_f\geq\frac{1}{134K^{4n-1}}.
$$
\end{Thm}


Our main aim  of this article is to consider the following natural questions.

\bprob\label{HLP-prob1}
Can we establish some new versions of Landau-type theorems for pluriharmonic mapping,
which are different with Theorems~D and F? 
\eprob

\bprob\label{HLP-prob2}
Can we improve Theorems~E and G? 
\eprob

The paper is organized as follows.  In Section \ref{HLP-sec2}, we present the main results of this paper. In Theorems \ref{HLP-th1}-\ref{HLP-th5}, we present an affirmative answer to Problem \ref{HLP-prob1}.  
In Section \ref{HLP-sec4}, we state and prove four related theorems which improve two results of Chen and Gauthier \cite{CHH2011}.
In particular, Theorems \ref{HLP-th4} and \ref{HLP-th6} provide an affirmative answer to Problem \ref{HLP-prob2}.

\setcounter{equation}{0}
\section{Main Results}\label{HLP-sec2}

We first state two Landau-type theorems of pluriharmonic mapping with bounded distortion, which are analogues versions of Theorem~D. 

\bthm\label{HLP-th1}
Let $f\in {\mathcal{PH}}^{1}(B^{n}) $ such that $f(0) = 0$ and $\Lambda_f(z)\leq\Lambda$ for $z\in B^{n}$.
Then $\Lambda\geq 1$ and $f\in {\mathcal S}_{PH}(\rho_1;R_1)$, where
$$ \rho_1=\frac{\pi}{4\Lambda_f^{2n-1}(0)(\Lambda_f(0)+\Lambda)} ~\mbox{ and }~R_1=\frac{\pi}{8\Lambda_f^{4n-2}(0)(\Lambda_f(0)+\Lambda)}.
$$
If, in addition, $\Lambda_f(0)=1$, then $f\in {\mathcal S}_{PH}(\rho_1';R_1')$, where
$$
\rho_1'=\frac{\pi}{4(1+\Lambda)}~\mbox{ and }~ R_1'=\frac{\pi}{8(1+\Lambda)}.
$$
\ethm



Next, we state a Landau-type theorem of locally univalent pluriharmonic $K$-mapping with bounded distortion, which is the analogues version of Theorem~F. 

\bthm\label{HLP-th5}
Let $f\in {\mathcal{PH}}_{loc,\,K}^{1}(B^{n})$, $n > 1$, such that $f(0) = 0$ and $\Lambda_f(z)\leq\Lambda$ for $z\in B^{n}$.
Then $\Lambda\geq 1/K^{2n-1}$ and that  there exists a domain $\Omega\subset B^n(0, \rho_2)$ such that $0\in\Omega$
and $f$ maps $\Omega$ onto a ball $B^n(0,R_2)$ injectively, where
$$
\rho_2 = \frac{\pi}{4 K^{2n-1}(K+\Lambda)}~\mbox{ and }~ R_2 = \frac{\pi}{8 K^{4n-2}(K+\Lambda)} .
$$
\ethm

Now, we state a Bloch-type theorem of pluriharmonic $K$-quasiregular mapping. 

\bthm\label{HLP-th4}
Let $f\in {\mathcal{PH}}_K(B^{n})$ and $n > 1$. Then $f(B^{n})$ contains a disk of radius $b_f$, with
\begin{eqnarray*}
b_f\geq R_3 =\frac{3-2\sqrt{2}}{8}\left (\frac{\pi}{K^{4n-1}}\right).
\end{eqnarray*}
\ethm

\br
Note that $R_3$ in Theorem \ref{HLP-th4} has a simple expression, and
$$ \frac{3-2\sqrt{2}}{8} > 0.0214466
$$
which is independent of $n$. Moreover, from Theorem~E, 
we have
\begin{eqnarray*}
R_n' &=& \frac{\pi}{8me(2n-1)}\Big(\frac{\pi}{4}\Big)^{4n-1}\frac{1}{K^{4n-1}}(1+o(1))\\
&\approx &\frac{1}{33.6 (2n-1)e}\Big(\frac{\pi}{4}\Big)^{4n-1}\frac{\pi}{K^{4n-1}}(1+o(1)),
\end{eqnarray*}
the number $\frac{1}{33.6 (2n-1)e}\Big(\frac{\pi}{4}\Big)^{4n-1}$ is a decreasing function of $n>1$. In particular, when $n=2$, we have
\begin{eqnarray*}
\frac{1}{33.6 (2n-1)e}\Big(\frac{\pi}{4}\Big)^{4n-1}\approx 0.0006727816777,
\end{eqnarray*}
which implies that $R_3>31.8\cdot R_n'$ for all $n>1$ which shows that the constant in Theorem \ref{HLP-th4} has a more substantial improvement when 
compared to the estimate of Chen and Gauthier \cite{CHH2011}, namely, Theorem~E. 
\er

\bcor\label{HLP-cor3}
Let $f\in {\mathcal{PH}}_K^{1}(B^{n})$, $n > 1$ and such that $\Lambda_f(0) = 1$. Then $f(B^{n})$ contains a disk of radius $b_f$, with
\begin{eqnarray*}
b_f\geq R_3' =\frac{(2-\sqrt{2})\pi}{8(1+(\sqrt{2}+1)K)}\geq\frac{3-2\sqrt{2}}{8}\cdot\frac{\pi}{K}>0.0214466\cdot\frac{\pi}{K}.
\end{eqnarray*}
\ecor

Finally, we state a Bloch-type theorem of locally univalent pluriharmonic $K$-mappings, which also has bigger improvement compared to the estimate of Theorem~G 
due to Chen and Gauthier \cite{CHH2011}.

\bthm\label{HLP-th6}
Let $f$ be a locally univalent pluriharmonic $K$-mapping of the unit
ball $B^{n}$ into $\mathbb{C}^{n}$, $n > 1$, such that ${\rm det} J_f (0) =1$.  Then $f(B^{n})$ contains a disk of radius $b_f$, with
\begin{eqnarray*}
b_f\geq \frac{3-2\sqrt{2}}{8}\cdot\frac{\pi}{K^{4n-1}}>\frac{1}{14.8 K^{4n-1}}>9\cdot\frac{1}{134 K^{4n-1}}.
\end{eqnarray*}
\ethm


\section{Proofs of the main results}\label{HLP-sec4}

First we recall the following well-known lemma. 

\begin{Lem}\label{HLP-lem1} {\rm (cf. Chen et al. \cite[Theorem 4]{CHH2011})} 
Let $f$ be a pluriharmonic mapping of $B^{n}$ into $B^{m}$. Then
\begin{equation}\label{eqn31}
\Lambda_f(z)\leq\frac{4}{\pi}\frac{1}{1-|z|^2} ~\mbox{ for }~  z\in  B^{n}.
\end{equation}
If $f(0)=0$, then
\begin{equation}\label{eqn32}
|f(z)|\leq\frac{4}{\pi}\arctan|z| \leq \frac{4}{\pi}|z|
~\mbox{ for }~  z\in B^{n}.
\end{equation}
\end{Lem}

\subsection{Proof of Theorem \ref{HLP-th1}}
By assumption $f(0) = 0$, $\Lambda_f(z)\leq\Lambda$ for $z\in B^{n}$ and ${\rm det} J_f (0) =1$. Thus,
we have 
$$
\Lambda\geq\Lambda_f(0)\geq\lambda_f(0)\geq \frac{|{\rm det} J_f (0)|}{\Lambda_f^{2n-1}(0)}=\frac{1}{\Lambda_f^{2n-1}(0)},
$$
which implies that  $\Lambda\geq\Lambda_f(0)\geq 1.$

Fix two distinct points $z', z{''}\in B^n (0, \rho_1)$, let $z{''}-z'=|z{''}-z'|\theta$ and define the pluriharmonic mapping $\phi_{\theta}$ by
\begin{equation*}
\phi_{\theta}(z)=(f_z(z)-f_z(0))\theta+(f_{\overline{z}}(z)-f_{\overline{z}}(0))\overline{\theta}.
\end{equation*}
Then, it follows from the definition of $\Lambda_f(z)$ that
\begin{eqnarray*}
|\phi_{\theta}(z)|&\leq &\Lambda_f(0)+\Lambda_f(z)\leq \Lambda_f(0)+\Lambda  ~\mbox{ for }~  z\in B^n.
\end{eqnarray*}
Note that $\phi_{\theta}(0)=0$. Using \eqref{eqn32}, we obtain
\begin{eqnarray*}
|\phi_{\theta}(z)|\leq\frac{4(\Lambda_f(0)+\Lambda)|z|}{\pi} ~\mbox{ for }~ z\in B^n.
\end{eqnarray*}
Consequently,
\begin{eqnarray*}
|f(z'')-f(z')|&=&\left |\int_{[z', z'']}f_z(z)\,dz+f_{\overline{z}}(z)\,d\overline{z}\right |\nonumber\\
&\geq&\left |\int_{[z', z'']}f_z(0)\,dz+f_{\overline{z}}(0)\,d\overline{z}\right |\\
&&-\left |\int_{[z', z'']}(f_z(z)-f_z(0))\,dz+(f_{\overline{z}}(z)-f_{\overline{z}}(0))\,d\overline{z}\right|,\nonumber\\
&\geq & |z'-z''|\lambda_f(0)-\int_{[z', z'']}|\phi_{\theta}(z)|\,ds\nonumber\\
&>& \frac{|z'-z''|}{\Lambda_f^{2n-1}(0)}-\frac{4(\Lambda_f(0)+\Lambda)\rho_1}{\pi}|z'-z''|\\
&=& \frac{4(\Lambda_f(0)+\Lambda)}{\pi}|z'-z''|\left [ \frac{\pi}{4\Lambda_f^{2n-1}(0)(\Lambda_f(0)+\Lambda)}-\rho_1\right ]=0,
\end{eqnarray*}
which implies that $f(z')\neq f(z'')$ and thus,  $f$ is univalent in $B^n (0, \rho_1)$.

Now let $z'\in \partial B^n (0, \rho_1)$. As $f(0)=0$, as above, we have
\begin{eqnarray*}
|f(z')| 
&\geq&\left |\int_{[0, z']}f_z(0)\,dz+f_{\overline{z}}(0)\,d\overline{z}\right |-\left |\int_{[0, z']}(f_z(z)-f_z(0))\,dz+(f_{\overline{z}}(z)-f_{\overline{z}}(0))\,d\overline{z}\right|\nonumber\\
&\geq& \lambda_f(0)\rho_1-\int_0^{\rho_1}\frac{4(\Lambda_f(0)+\Lambda) r}{\pi}\,dr\nonumber\\
&\geq& \frac{\rho_1}{\Lambda_f^{2n-1}(0)}-\frac{2(\Lambda_f(0)+\Lambda)\rho_1^2}{\pi}= \frac{\pi}{8\Lambda_f^{4n-2}(0)(\Lambda_f(0)+\Lambda)}=R_1,
\end{eqnarray*}
which shows that $f(B^n (0, \rho_1))$ covers the ball $B^n (0, R_1)$, and the proof is complete. \hfill $\Box$

\subsection{Proof of Theorem \ref{HLP-th5}}
Let $f\in {\mathcal{PH}}_{loc,\,K}^{1}(B^{n})$ with $n > 1$, i.e. $f$ is a locally univalent pluriharmonic $K$-mapping of the unit
ball $B^{n}$ into $\mathbb{C}^{n}$ such that ${\rm det} J_f (0)=1$. As $\Lambda_f(z)\leq\Lambda$ for $z\in B^{n}$ and
$$\lambda_f(0)\geq \frac{|{\rm det} J_f (0)|}{K^{2n-1}}=\frac{1}{K^{2n-1}},
$$
it follows that
$$
\Lambda_f(0)=|J_f (0)| \leq K|{\rm det} J_f (0)|^{1/(2n)}=K
~\mbox{ and }
\Lambda\geq\Lambda(0)\geq\lambda_f(0)\geq \frac{1}{K^{2n-1}}.
$$

As with the previous theorem, let  $z', z{''}\in B^n (0, \rho_2)$ be two distinct points,  $z{''}-z'=|z{''}-z'|\theta$ and
\begin{equation*}
\phi_{\theta}(z)=(f_z(z)-f_z(0))\theta+(f_{\overline{z}}(z)-f_{\overline{z}}(0))\overline{\theta}.
\end{equation*}
Then,  the definition of $\Lambda_f(z)$ gives that
\begin{eqnarray*}
|\phi_{\theta}(z)|&\leq &\Lambda_f(0)+\Lambda_f(z)\leq K+\Lambda ~\mbox{ for }~  z\in B^n.
\end{eqnarray*}
Again, as $\phi_{\theta}(0)=0$, \eqref{eqn32} implies that
\begin{eqnarray*}
|\phi_{\theta}(z)|\leq\frac{4(K+\Lambda) |z|}{\pi} ~\mbox{ for }~ z\in B^n.
\end{eqnarray*}
Then, using the analogous proof of Theorem \ref{HLP-th1}, we have 
\begin{eqnarray*}
|f(z'')-f(z')|&\geq& |z'-z''|\lambda_f(0)-\int_{[z', z'']}|\phi_{\theta}(z)|\,ds\nonumber\\
&>& \frac{|z'-z''|}{K^{2n-1}}-\frac{4(K+\Lambda)\rho_2}{\pi}|z'-z''|=0,
\end{eqnarray*}
which shows that $f$ is univalent in $B^n (0, \rho_2)$.

Now let $z'\in \partial B^n (0, \rho_2)$ and observe that $f(0)=0$. As in the proof of Theorem \ref{HLP-th1}, we have
\begin{eqnarray*}
|f(z')|
&\geq& \lambda_f(0)\rho_2-\int_0^{\rho_2}\frac{4(K+\Lambda) r}{\pi}\,dr\nonumber\\
&=& \frac{\rho_2}{K^{2n-1}}-\frac{2(K+\Lambda)\rho_2^2}{\pi}= \frac{\pi}{8 K^{4n-2}(K+\Lambda)}=R_2.
\end{eqnarray*}
This shows that $f(B^n (0, \rho_2))$ covers the ball $B^n (0, R_2)$, and the proof is complete. \hfill $\Box$

\subsection{Proof of Theorem \ref{HLP-th4}}
By means of Theorem \ref{HLP-th1},  we may use arguments similar to those in the proof of \cite[Theorem 6]{CHH2011}. For the sake of readability,
we provide the details.

Let $f\in {\mathcal{PH}}_K(B^{n})$ and $n > 1$. Without loss of generality, we assume that $f$ is pluriharmonic on $\overline{B}^n$. Then $(1-|z|)^{2n}|{\rm det} J_f (z)|$ is
continuous and bounded  on $\overline{B}^n$. Moreover,
\begin{eqnarray*}
(1-|z|)^{2n}|{\rm det} J_f (z)|\Big|_{z=0}=|{\rm det} J_f (0)|=1~\mbox{ and }~\lim\limits_{|z|\to 1^-}(1-|z|)^{2n}|{\rm det} J_f (z)|=0.
\end{eqnarray*}
This implies that there exists a point $z'\in B^n$ such that
\begin{eqnarray}
(1-|z'|)^{2n}|{\rm det} J_f (z')|=1
\label{eqn33}
\end{eqnarray}
and
\begin{eqnarray*}
(1-|z|)^{2n}|{\rm det} J_f (z)|\leq 1~\mbox{ for }~ |z'|=r\leq |z|\leq 1.
\label{eqn34}
\end{eqnarray*}
In particular, we have
\begin{eqnarray}
|{\rm det} J_f (z)|\leq |{\rm det} J_f (z')|~\mbox{ for }~|z|=r.
\label{eqn35}
\end{eqnarray}

In the following, we need to consider two cases,

\vspace{8pt}
\noindent
{\bf Case 1.  $r>0$}: Fix a point $z_0$ with $0<|z_0|\leq r$ and consider $\Lambda_f(z_0)=|f_z(z_0)\theta+f_{\overline{z}}(z_0)\overline{\theta}|$ with $\theta\in\partial B^n$. Define the function $\phi$ by
\begin{eqnarray*}
\phi(\zeta)=f_z(\zeta z_0/|z_0|)\theta+f_{\overline{z}}(\zeta z_0/|z_0|)\overline{\theta}~\mbox{ for }~  \zeta \in \ID.
\end{eqnarray*}

Since $\phi$ is harmonic, by the maximum principle, there exists a point $\zeta'$ with $|\zeta'|=r$, such that 
\begin{eqnarray*}
\Lambda_f(z_0)&=&|\phi(|z_0|)|\leq |f_z(\zeta' z_0/|z_0|)\theta+f_{\overline{z}}(\zeta' z_0/|z_0|)\overline{\theta}|.
\end{eqnarray*}
Let $z''=\zeta' z_0/|z_0|$. Note that $|z''|=r$, by the definition of $K$-mappings and (\ref{eqn35}), we have
\begin{eqnarray*}
\Lambda_f(z_0)&\leq & |f_z(z'')\theta+f_{\overline{z}}(z'')\overline{\theta}|\\
&\leq & \Lambda_f(z'')=|J_f(z'')|\leq K|{\rm det} J_f(z'')|^{1/(2n)}\\
&\leq & K|{\rm det} J_f(z')|^{1/(2n)}.
\end{eqnarray*}
On the other hand, by the definition of $K$-mappings, (\ref{eqn33}) and $|{\rm det} J_f (0)|=1$, we have
\begin{eqnarray*}
\Lambda_f(0)&=&|J_f(0)|\leq K |{\rm det} J_f(0)|^{1/(2n)}=K\\
&=& K (1-r)|{\rm det} J_f(z')|^{1/(2n)}\\
&\leq & K|{\rm det} J_f(z')|^{1/(2n)}.
\end{eqnarray*}
Thus we conclude that
\begin{eqnarray}
\Lambda_f(z)\leq |{\rm det} J_f (z')|^{1/(2n)}~\mbox{ for }~ |z|\leq r.
\label{eqn36}
\end{eqnarray}

Now we define the functions $g$ and $F$ by
\begin{eqnarray*}
g(\zeta)=z'+(2-\sqrt{2})(1-r)\zeta ~\mbox{ and }~ F(\zeta)=\frac{1}{2-\sqrt{2}}(f(g(\zeta))-f(z'))
\end{eqnarray*}
for $\zeta\in B^n$. Then
\begin{eqnarray*}
F(0)=0 ~\mbox{ and }~ |{\rm det} J_F (0)|=(1-r)^{2n}|{\rm det} J_f (z')|=1.
\end{eqnarray*}

If $|g(\zeta)|\leq r$, by (\ref{eqn36}) and (\ref{eqn33}), we have
\begin{eqnarray*}
\Lambda_F(\zeta)&=& (1-r)\Lambda_f(g(\zeta)) 
\leq K(1-r)|{\rm det} J_f (z')|^{1/(2n)}=K;
\end{eqnarray*}
and if $|g(\zeta)|\geq r$,
\begin{eqnarray}
\Lambda_F(\zeta)&=& (1-r)\Lambda_f(g(\zeta))\leq  K(1-r)|{\rm det} J_f (g(\zeta))|^{1/(2n)}\nonumber\\
&=&K\left (\frac{1-r}{1-|g(\zeta)|}\right )(1-|g(\zeta)|)|{\rm det} J_f (g(\zeta))|^{1/(2n)}\nonumber\\
&\leq & K\left (\frac{1-r}{1-|g(\zeta)|}\right )\leq \frac{K(1-r)}{1-r-(2-\sqrt{2})(1-r)|\zeta|}\nonumber\\
& = &\frac{K}{1-(2-\sqrt{2})|\zeta|}. 
\label{eqn37}
\end{eqnarray}

\vspace{8pt}
\noindent
{\bf Case 2.  $r=0$}:  Consider the functions  $g$ and $F$ defined as above. Then $|g(\zeta)|\geq r=0$ and it follows from (\ref{eqn37}) that
\begin{eqnarray*}
\Lambda_F (\zeta)\leq \frac{K}{1-(2-\sqrt{2})|\zeta|} ~\mbox{ for }~ \zeta\in B^n.
\end{eqnarray*}
Thus we conclude that
\begin{eqnarray*}
\Lambda_F (\zeta)\leq \frac{K}{1-(2-\sqrt{2})|\zeta|}<(\sqrt{2}+1) K ~\mbox{ for }~ \zeta\in B^n.
\end{eqnarray*}
In particular, $\Lambda_F (0)\leq K$.

Now, applying Theorem \ref{HLP-th1} to the mapping $F(z)$, we see that $F(B^{n})$ contains a schlicht ball with center $0$ and radius 
\begin{equation*}
R=\frac{\pi}{8K^{4n-1}(K+(\sqrt{2}+1)K)}=\frac{1}{8(2+\sqrt{2})}\cdot\frac{\pi}{K^{4n-1}}.
\end{equation*}
Consequently, $f(B^{n})$ contains a schlicht ball with center $f(z')$ and radius $R_3$, that is,
\begin{equation*}
b_f\geq R_3=(2-\sqrt{2})R =\frac{3-2\sqrt{2}}{8}\cdot\frac{\pi}{K^{4n-1}}>0.0214466\cdot\frac{\pi}{K^{4n-1}}.
\end{equation*}
This proves the theorem. \hfill $\Box$

\subsection{Proof of Theorem \ref{HLP-th6}}
By means of Theorem \ref{HLP-th5} (replacing Theorem \ref{HLP-th1}), and the analogous proof of Theorem \ref{HLP-th4}, we may finish the proof  of Theorem \ref{HLP-th6}. \hfill $\Box$
\vskip 10mm

\subsection*{Acknowledgments}
This research of the first author was partly supported by Guangdong Natural Science Foundations (Grant No. 2021A1515010058).
The work of the second author was supported by Mathematical Research Impact Centric Support (MATRICS) of the Department of Science and Technology (DST), India  (MTR/2017/000367).


\end{document}